\documentclass[draft]{amsart}

\usepackage{amssymb}
\usepackage{url}
\usepackage{amscd}
\usepackage{amsthm}
\usepackage{amsmath}

\theoremstyle{plain}

\theoremstyle{remark}

\begin{document}

\date{}

\title[The Graham conjecture implies the Erd\"{o}s-Tur\'{a}n conjecture]
{The Graham conjecture implies the Erd\"{o}s-Tur\'{a}n conjecture}
\author{Liangpan Li}

\date{April 4, 2007}

\address{Department of Mathematics,
Shanghai Jiaotong University, Shanghai 200240, People's Republic of
China}

\email{liliangpan@sjtu.edu.cn}


\begin{abstract}
Erd\"{o}s and Tur\'{a}n once conjectured that any set
$A\subset\mathbb{N}$ with $\sum_{a\in A}{1}/{a}=\infty$ should
contain infinitely many progressions of arbitrary length $k\geq3$.
For the two-dimensional case Graham conjectured that if $B\subset
\mathbb{N}\times\mathbb{N}$ satisfies
$$\sum\limits_{(x,y)\in B}\frac{1}{x^2+y^2}=\infty,$$ then for any $s\geq2$, $B$
contains an $s\times s$ axes-parallel grid. In this paper it is
shown that if the Graham conjecture is true for some $s\geq2$, then
the Erd\"{o}s-Tur\'{a}n conjecture is true for $k=2s-1$.
\end{abstract}

\maketitle

\section{Introduction}
One famous conjecture of Erd\"{o}s and Tur\'{a}n \cite{erdos Turan}
asserts that any set $A\subset\mathbb{N}$ with $\sum_{a\in
A}{1}/{a}=\infty$ should contain infinitely many progressions of
arbitrary length $k\geq3$. There are two important progresses
towards this direction due to Szemer\'{e}di \cite{szemeredi} and
Green and Tao \cite{green tao} respectively, which assert that if
$A$ has positive upper density or $A$ is the set of the prime
numbers, then $A$ contains infinitely many progressions of arbitrary
length.

If one considers the similar question in the two-dimensional plane,
Graham \cite{graham} conjectured that if $B\subset
\mathbb{N}\times\mathbb{N}$ satisfies
$$\sum\limits_{(x,y)\in B}\frac{1}{x^2+y^2}=\infty,$$ then $B$
contains the four vertices of an axes-parallel square. More
generally, for any $s\geq2$ it should be true that $B$ contains an
$s\times s$ axes-parallel grid. Furstenberg and Katznelson
\cite{Furstenberg Katznelson} proved the two-dimensional
Szemer\'{e}di theorem, that is, any set $B\subset
\mathbb{N}\times\mathbb{N}$ with positive upper density contains an
$s\times s$ axes-parallel grid. In another words, such a set $B$
contains any finite pattern.

The purpose of this paper is to show that if the Graham conjecture
is true, then the Erd\"{o}s-Tur\'{a}n conjecture is also true.

\section{The Graham conjecture implies the Erd\"{o}s-Tur\'{a}n conjecture}

Suppose that the Erd\"{o}s-Tur\'{a}n conjecture is false for $k=3$.
Then there exists a set
$$A=\{a_{1}<a_{2}<a_{3}<\cdots\} \subset\mathbb{N}$$ with
$\sum_{n\in\mathbb{N}}{1}/{a_{n}}=\infty$ such that $A$ contains no
arithmetic progression of length 3. Define a set
$B\subset\mathbb{N}\times\mathbb{N}$ by
$$B=\Big\{(a_{n}+m,m):n\in\mathbb{N},m\in\mathbb{N}\Big\}.$$
Then \begin{align*}\sum_{(x,y)\in B}\frac{1}{x^{2}+y^{2}}
&=\sum_{n\in\mathbb{N}}\sum_{m\in\mathbb{N}}\frac{1}{(a_{n}+m)^{2}+m^{2}}\\
&\geq\sum_{n\in\mathbb{N}}\sum_{m=1}^{a_{n}}\frac{1}{(a_{n}+m)^{2}+m^{2}}\\
&\geq\sum_{n\in\mathbb{N}}\frac{a_{n}}{(a_{n}+a_{n})^{2}+a_{n}^{2}}\\
&=\sum_{n\in\mathbb{N}}\frac{1}{5a_{n}}\\&=\infty.\end{align*}

In the sequel we indicate that $B$ contains no square and argue it
by contradiction. This would mean that the Graham conjecture is
false for $s=2$. Suppose that for some $n,m,l\in\mathbb{N}$, $B$
contains a square of the following form:
\begin{align*}
(a_{n}+m,m+l),\ \  &(a_{n}+m+l,m+l),\\
(a_{n}+m,m),\,   \qquad  &(a_{n}+m+l,m).\end{align*}  It follows
easily from the construction of $B$ that $a_{n}-l,a_{n},a_{n}+l\in
A$, which yields a contradiction since $A$ contains no arithmetic
progression of length 3 according to the initial assumption.

Similarly, if the Graham conjecture is true for some $s\geq2$, then
the Erd\"{o}s-Tur\'{a}n conjecture is true for $k=2s-1$. The
interested reader can easily provide a proof.

\section{Concluding Remarks}
Let $r(k,N)$ be the maximal cardinality of a subset $A$ of
$\{1,2,\ldots,N\}$ which is free of $k$-term arithmetic
progressions. Behrend \cite{behrend} and Rankin \cite{Rankin} had
shown that
$$r(k,N)\geq N\cdot\exp(-c(\log N)^{1/(k-1)}).$$
Similarly, let $\widetilde{r}(s,N)$ be the maximal cardinality of a
subset $B$ of $\{1,2,\ldots,N\}^{2}$ which is free of $s\times s$
axes-parallel grids. For any set $A\subset\{1,2,\ldots,N\}$, define
$$\Theta(A)=\{(a+m,m):a\in A,m=1,2,\ldots,N\}\subset\{1,2,\ldots,2N\}^{2}.$$
Following the discussion in Section 2, one can easily deduce that if
$A$ is free of $2s-1$ term of arithmetic progression, then
$\Theta(A)$ is free of $s\times s$ axes-parallel grid. Hence
\begin{align*}{\widetilde{r}(s,2N)}&\geq{r(2s-1,N)}\cdot{N}\\
&\geq N^{2}\exp(-c(\log N)^{1/(2s-2)}).\end{align*}

We end this paper with a question. Does the Erd\"{o}s-Tur\'{a}n
conjecture imply the Graham conjecture?

\end{document}